\documentclass[11pt]{article}
\usepackage{graphicx}
\usepackage{amsmath,amsthm,amssymb}

\numberwithin{equation}{section}

\begin{document}


\title{\bf Space-Time Foam Differential Algebras of Generalized Functions and a Global
Cauchy-Kovalevskaia Theorem}%

\author{Elem\'{e}r E Rosinger}

\maketitle

\begin{center}
Department of Mathematics and Applied Mathematics,\\
University of Pretoria,\\
0002, Pretoria, South Africa
\\
\verb"eerosinger@hotmail.com"
\end{center}

\begin{abstract}

The new {\it global} version of the Cauchy-Kovalevskaia theorem presented here is a
strengthening and extension of the {\it regularity} of similar global solutions obtained
earlier by the author. Recently the space-time foam differential algebras of generalized
functions with \textit{dense} singularities were introduced. A main motivation for these
algebras comes from the so called space-time foam structures in General Relativity, where the
set of singularities can be dense. A variety of applications of these algebras have been
presented elsewhere, including in de Rham cohomology, Abstract Differential Geometry, Quantum
Gravity, etc. Here a global Cauchy-Kovalevskaia theorem is presented for arbitrary
analytic nonlinear systems of PDEs. The respective global generalized solutions are analytic
on the whole of the domain of the equations considered, except for singularity sets which are
closed and nowhere dense, and which upon convenience can be chosen to have zero Lebesgue
measure. \\
In view of the severe limitations due to the polynomial type growth conditions in the
definition of Colombeau algebras, the class of singularities such algebras can deal with is
considerably limited. Consequently, in such algebras one cannot even formulate, let alone
obtain, the global version of the Cauchy-Kovalevskaia theorem presented in this paper.

\bigskip

\bigskip

\begin{quote}
        ``We do not possess any method at all to derive systematically
        solutions that are free of singularities...''

        \medskip

        Albert Einstein
        \\
        \textit{The Meaning of Relativity}
        \\
        Princeton Univ. Press, 1956, p. 165
\end{quote}

\end{abstract}

\section{Algebras of Generalized functions with Dense Singularities, or Space-Time Foam
Algebras}

\subsection{Families of Dense Singularities in Euclidean Spaces}

In this paper, following Rosinger [9-11,13,15], we consider differential algebras of
generalized functions - called {\it space-time foam} algebras - which have significantly
strengthened and extended properties with respect to the {\it singularities} they can deal
with. Namely, this time the singularities can be arbitrary, including \textit{dense} sets, and
the only condition they have to satisfy is that their complementaries, that is, the set of
nonsingular points, be also dense. This, among others, allows for singularity sets with a
cardinal \textit{larger} than that of the set of nonsingular points. For instance, in the case
the domain is the real line, the set of singularities can be given by the uncountable set of
all the irrational numbers, since its complementary, the set of rational numbers, is still
dense, although it is only countable. \\

These space-time foam algebras are instances of the earlier nonlinear algebraic theory of
generalized functions introduced and developed in Rosinger [1-8,13-15], Rosinger \& Walus
[1,2], Mallios \& Rosinger [1], Mallios [1], see 46F30 at www.ams.org/msc/46Fxx.html \\
This general nonlinear algebraic theory has so far exhibited as particular cases a number of
differential algebras of generalized functions, among them, the Colombeau algebras, see
Grosser et.al. [p. 7]. \\

The space-time foam algebras in this paper are able to deal with by far the {\it largest}
class of singularities so far in the literature. \\
This fact proves to have useful {\it existence} and {\it regularity} consequences when it
comes to the solutions in the global version of the Cauchy-Kovalevskaia theorem. \\

On the other hand, in view of the severe limitations due to the polynomial type growth
conditions in the definition of Colombeau algebras, the class of singularities such algebras
can deal with is considerbaly limited. Consequently, in such algebras one cannot even
formulate, let alone obtain, the global version of the Cauchy-Kovalevskaia theorem presented
in this paper. \\

In this section, following Rosinger [9-11] where they were first introduced, we recall in
short the construction of these new, namely, space-time foam algebras of generalized functions.
For that purpose, first we have to introduce the {\it families of singularities} such algebras
can deal with. \\

Let our underlying topological space $X$ be any nonvoid open subset of $\mathbb{R}^n$. The
general case in the construction of space-time foam algebras, namely when $X$ is any finite
dimensional smooth manifold, is presented in Rosinger [11], and rather surprisingly, it does
not lead to any additional technical difficulties. This fact, in addition to the far larger
class of singularities it can handle, is one of the advantages of the space-time foam algebras
when compared, for instance, with the Colombeau algebras. \\

We shall consider various families of singularities in $X$, each such family being given by a
corresponding set $\mathcal{S}$ of subsets $\Sigma \subset X$, with each such subset $\Sigma$
describing a possible set of singularities of a certain given generalized function, or in
particular, generalized solution.

\medskip
The \textit{largest} family of singularities $\Sigma \subset X$ which we can consider is given
by

\begin{equation}\label{label1.1}
    \mathcal{S}_{\mathcal{D}} ~=~ \{~ \Sigma \subset X ~~|~~
                            X\setminus \Sigma \textnormal{ is dense in } X ~\}
\end{equation}

In this way, the various families $\mathcal{S}$ of singularities $\Sigma \subset X$ which we
shall deal with, will each satisfy the condition $\mathcal{S}\subseteq \mathcal{S}_
{\mathcal{D}}$.

\bigskip
Regarding the treatment of {\it singularities} in a Differential Geometric context, it should
be noted that a major interest in large, possibly \textit{dense} sets of singularities comes
from general relativity, see Finkelstein, Geroch [1,2], Heller [1-3], Heller \& Sasin [1-3],
Heller \& Multarzynski \& Sasin, Gurszczak \& Heller, or Mallios [1-6], Mallios \& Rosinger
[1-3]. \\

In this respect we note that, according to the strongest earlier corresponding result, see
Heller [2], Heller \& Sasin [2], the family of singularities $\mathcal{S}$ could only be
composed from one single closed nowhere dense $\Sigma$, which in addition, had to be in the
boundary of $X$. \\

On the other hand, in Mallios \& Rosinger [1] - which except for Mallios \& Rosinger [2,3],
did treat the most general type of singularities - the family $\mathcal{S}$ could already
contain {\it all closed and nowhere dense} subsets $\Sigma$ in $X$. And then finally, in
Mallios \& Rosinger [2,3] the {\it largest} class of singularities so far, namely such as in
this paper, thus in particular, dense singularities as well, are treated.

\bigskip

For earlier developments regarding the possible treatment of singularities in a Differential
Geometric context one can consult, for instance, Sikorski, Kirillov [1,2], Mostow, or Souriau
[1,2].
And it should be mentioned that, as seen in Finkelstein, and especially in Geroch [1,2], the
issue of singularities has for a longer time been of fundamental importance in General
Relativity.

\bigskip

In this paper, as in Rosinger [9-11,13,15] and Mallios \& Rosinger [2,3], the family
$\mathcal{S}$ of singularities can be any subset of $\mathcal{S}_{\mathcal{D}}$ in
(\ref{label1.1}). Among other ones, two such families which will be of interest are the
following

\begin{equation}\label{label1.2}
    \mathcal{S}_{nd} ~=~ \{~ \Sigma \subset X ~~|~~
                \Sigma \textnormal{ is closed and nowhere dense in } X ~\}
\end{equation}

and

\begin{equation}\label{label1.3}
    \mathcal{S}_{\textit{Baire I}} ~=~ \{~ \Sigma \subset X ~~|~~
                   \Sigma \textnormal{ is of first Baire category in } X ~\}
\end{equation}

Obviously

\begin{equation}\label{label1.4}
    \mathcal{S}_{nd} ~\subset~ \mathcal{S}_{\textit{Baire I}} ~\subset~
    \mathcal{S}_{\mathcal{D}}
\end{equation}

\subsection{Asymptotically Vanishing Ideals}

Let us now for convenience recall shortly the idea of the construction introduced in Rosinger
[9-11]. There are \textit{two} basic ingredients involved. First, we take any family
$\mathcal{S}$ of singularity sets $\Sigma \subset X$, family which satisfies the conditions

\begin{equation}\label{label1.5}
\begin{split}
    & \forall ~ \Sigma \in \mathcal{S}\, : \\&
    \quad X\setminus \Sigma \textnormal{ is dense in } X
\end{split}
\end{equation}

and

\begin{equation}\label{label1.6}
\begin{split}
        & \forall ~ \Sigma , \Sigma \, '\in \mathcal{S}\, :\\&
        \exists ~ \Sigma \, '' \in \mathcal{S}\, :\\&
        \quad \Sigma \cup \Sigma \, ' \subseteq \Sigma \, ''
\end{split}
\end{equation}

Clearly, both families $\mathcal{S}_{nd}$ and $\mathcal{S}_{\textit{Baire I}}$   satisfy
conditions (\ref{label1.5}) and  (\ref{label1.6}).\\

Now, as the second ingredient, and so far independently of $\mathcal{S}$ above, we take any
right directed partial order $L=(\Lambda, \leq)$. In other words, $L$ is such that for each
$\lambda, \lambda \, ' \in \Lambda$, there exists $\lambda \, '' \in \Lambda$, with $\lambda,
\lambda\, ' \leq \lambda \, ''$.  The role of $L$ will become clear later, see for instance
Example 1 in section 2.\\

Although we shall only be interested in singularity sets $\Sigma \in
\mathcal{S}_{\mathcal{D}}$, the following ideal can be defined for any $\Sigma \subseteq X$.
Indeed, let us denote by

\begin{equation}\label{label1.7}
    \mathcal{J}_{L,\, \Sigma}(X)
\end{equation}

the \textit{ideal} in $(\mathcal{C}^{\infty}(X))^{\Lambda}$ of all the sequences of smooth
functions indexed by $\lambda \in \Lambda$, namely, $w=(w_{\lambda} \, | \, \lambda \in
\Lambda) \in (\mathcal{C}^{\infty}(X))^{\Lambda}$, sequences which \textit{outside} of the
singularity set $\Sigma$ will satisfy the \textit{asymptotic vanishing} condition

\begin{equation}\label{label1.8}
\begin{split}
        & \forall ~ x \in X \setminus \Sigma :\\ &
        \exists ~ \lambda \in \Lambda :\\&
        \forall ~ \mu \in \Lambda, \, \mu \geq \lambda :\\&
        \forall ~ p \in \mathbb{N}^n :\\&
        \quad D^p w_{\mu}(x)=0
\end{split}
\end{equation}

This means that the sequences of smooth functions $w=(w_{\lambda}\, | \, \lambda \in \Lambda)$
in the ideal $\mathcal{J}_{L,\, \Sigma}(X)$ may in a way \textit{cover} with their support the
singularity set $\Sigma$, and at the same time, they {\it vanish asymptotically} outside of it,
together with all their partial derivatives.
\\
In this way, the ideal $\mathcal{J}_{L,\, \Sigma}(X)$ carries in an \textit{algebraic} manner
the information on the singularity set $\Sigma$. Therefore, a \textit{quotient} in which the
factorization is made with such ideals may in certain ways \textit{do away with} singularities,
and do so through purely {\it algebraic} means, see (\ref{label1.11}), (\ref{label1.12})
below.

\bigskip

We note that the assumption about $L =(\Lambda, \leq)$ being right directed is used in proving
that $\mathcal{J}_{L,\, \Sigma}(X)$ is indeed an ideal, more precisely that, for $w, w\, ' \in
\mathcal{J}_{L,\, \Sigma}(X)$, we have $w + w\, ' \in \mathcal{J}_{L,\, \Sigma}(X)$. \\

Now, it is easy to see that for $\Sigma, \Sigma \, ' \subseteq X$, we have

\begin{equation}\label{label1.9}
    \Sigma ~\subseteq~ \Sigma \, ' ~~~\Longrightarrow~~~
              \mathcal{J}_{L,\, \Sigma}(X) ~\subseteq~ \mathcal{J}_{L,\, \Sigma \,'}(X)
\end{equation}

in this way, in view of (\ref{label1.6}), it follows that

\begin{equation}\label{label1.10}
    \mathcal{J}_{L,\, \mathcal{S}}(S) ~=~
                      \bigcup_{\Sigma \in \mathcal{S}}\mathcal{J}_{L,\, \Sigma}(X)
\end{equation}

is also an ideal in $(\mathcal{C}^{\infty}(X))^{\Lambda}$.

\subsection{Foam Algebras}

In view of the above, for $\Sigma \subseteq X$, we can define the algebra

\begin{equation}\label{label1.11}
    B_{L,\, \Sigma}(X) ~=~ (\mathcal{C}^{\infty}(X))^{\Lambda}/\mathcal{J}_{L,\, \Sigma}(X)
\end{equation}

However, we shall only be interested in singularity sets $\Sigma \in
\mathcal{S}_{\mathcal{D}}$, that is, for which $X\setminus \Sigma$ is {\it dense} in $X$. And
in such a case the corresponding algebra $B_{L,\, \Sigma}(X)$ will be called a \textit{foam
algebra.}

\subsection{ Multi-Foam Algebras}

With the given family $\mathcal{S}$ of singularities, and based on (\ref{label1.10}), we can
now associate the \textit{multi-foam algebra}

\begin{equation}\label{label1.12}
    B_{L,\, \mathcal{S}}(X) ~=~ (\mathcal{C}^{\infty}(X))^{\Lambda}/\mathcal{J}_{L,\,
                 \mathcal{S}}(X)
\end{equation}

\subsection{Space-Time Foam Algebras}

The foam algebras and the multi-foam algebras introduced above will for the sake of simplicity
be called together \textit{space-time foam algebras}. \\

Clearly, if the family $\mathcal{S}$ of singularities consists of one single singularity set
$\Sigma \in \mathcal{S}_{\mathcal{D}}$, that is, $\mathcal{S} = \{\, \Sigma \,\}$, then
conditions (\ref{label1.5}), (\ref{label1.6}) are satisfied, and in this particular case the
concepts of foam and multi-foam algebras are identical, in other words, $B_{L,\, \{\, \Sigma
\,\}}(X) = B_{L,\, Sigma}(X)$. This means that the concept of multi-foam algebra is more
general than that of foam algebra.

\bigskip
It is clear from their quotient construction that the space-time foam algebras are associative
and commutative. However, the above constructions can easily be extended to the case when,
instead of real valued smooth functions, we use smooth functions with values in an arbitrary
\textit{normed algebra}. In such a case the resulting space-time foam algebras will still be
associative, but in general they may be noncommutative.

\subsection{Space-Time Foam Algebras as Algebras of Generalized Functions}

The reason why we restrict ourself to singularity sets $\Sigma \in \mathcal{S}_{\mathcal{D}}$,
that is, to subsets $\Sigma \subset X$ for which $X\setminus \Sigma$ is dense in $X$, is due
to the implication, see further details in Rosinger [15], and for a full argument Rosinger [4,
chap. 3, pp. 65-119]

\begin{equation}\label{label1.13}
    X \setminus \Sigma \textnormal{ is dense in } X ~~~\Longrightarrow~~~
                  \mathcal{J}_{L,\, \Sigma}(X)
    \cap \mathcal{U}_{\Lambda}^{\infty}(X) ~=~ \{\, 0\,\}
\end{equation}

where $\mathcal{U}_{\Lambda}^{\infty}(X)$ denotes the \textit{diagonal} of the power
$(\mathcal{C}^\infty(X))^{\Lambda}$, namely, it is the set of all $u(\psi)=(\psi_{\lambda}\,
| \, \lambda \in \Lambda)$, where $\psi_{\lambda}=\psi$, for $\lambda\in \Lambda$, while
$\psi$ ranges over $\mathcal{C}^{\infty}(X)$. In this way, we have the algebra isomorphism
$\mathcal{C}^{\infty}(X) \ni \psi ~\longmapsto~ u(\psi) \in \mathcal{U}_{\Lambda}^{\infty}
(X)$.

\bigskip
This implication (\ref{label1.13}) follows immediately from the asymptotic vanishing condition
(\ref{label1.8}). Indeed, if $\psi \in {\cal C}^\infty ( X )$ and $u ( \psi ) \in
\mathcal{J}_{L,\, \Sigma}(X)$, then (1.8) implies that $\psi = 0$ on $X \setminus \Sigma$,
thus we must have $\psi = 0$ on $X$, since $X \setminus \Sigma$ was assumed to be dense in $X$.
It follows, therefore, that the ideal $\mathcal{J}_{L,\, \Sigma}(X)$ is \textit{off
diagonal}. \\

The importance of (\ref{label1.13}) is that, for $\Sigma \in \mathcal{S}_{\mathcal{D}}$, it
gives the following \textit{algebra embedding} of the smooth functions into foam algebras

\begin{equation}\label{label1.14}
    \mathcal{C}^{\infty}(X) \ni \psi ~\longmapsto~ u(\psi) +
                   \mathcal{J}_{L,\, \Sigma}(X) \in B_{L,\, \Sigma}(X)
\end{equation}

Now in view of (\ref{label1.10}), it is easy to see that (\ref{label1.13}) will as well yield
the \textit{off diagonality} property

\begin{equation}\label{label1.15}
    \mathcal{J}_{L,\, \mathcal{S}}(X) \cap \mathcal{U}_{\Lambda}^{\infty}(X) ~=~ \{\, 0 \,\}
\end{equation}

and thus similar with (\ref{label1.14}), we obtain the \textit{algebra embedding} of smooth
functions into multi-foam algebras

\begin{equation}\label{label1.16}
    \mathcal{C}^{\infty}(X) \ni \psi ~\longmapsto~ u(\psi) + \mathcal{J}_{L,\, \mathcal{S}}(X) \in B_{L,\, \mathcal{S}}(X)
\end{equation}

The algebra embeddings (\ref{label1.14}), (\ref{label1.16}) mean that the foam and multi-foam
algebras are in fact \textit{algebras of generalized functions}. Also they mean that the foam
and multi-foam algebras are unital, with the respective unit elements $u(1) +
\mathcal{J}_{L,\, \Sigma}(X), \, \, u(1) + \mathcal{J}_{L,\, \mathcal{S}}(X)$.

\bigskip
Further, the asymptotic vanishing condition (\ref{label1.8}) also implies quite obviously that,
for $\Sigma \subseteq X$, we have

\begin{equation}\label{label1.17}
    D^p \, \mathcal{J}_{L,\, \Sigma}(X) ~\subseteq~
                  \mathcal{J}_{L,\, \Sigma}(X), \textnormal{ for } p\in \mathbb{N}^n
\end{equation}

where $D^p$ denotes the termwise $p$-th order partial derivation of sequences of smooth
functions, applied to each such sequence in the ideal $\mathcal{J}_{L,\, \Sigma}(X)$. \\

Then again, in view of (\ref{label1.10}), we obtain

\begin{equation}\label{label1.18}
    D^p \, \mathcal{J}_{L,\, \mathcal{S}}(X) ~\subseteq~
       \mathcal{J}_{L,\, \mathcal{S}}(X), \textnormal{ for } p \in \mathbb{N}^n
\end{equation}

Now (\ref{label1.17}), (\ref{label1.18}) mean that the the foam and multi-foam algebras are in
fact \textit{differential algebras}, namely

\begin{equation}\label{label1.19}
    D^p \, B_{L,\, \Sigma}(X) ~\subseteq~
                            B_{L,\, \Sigma}(X), \textnormal{ for } p\in \mathbb{N}^n
\end{equation}

where $\Sigma \in \mathcal{S}_{\mathcal{D}}$, and furthermore we also have

\begin{equation}\label{label1.20}
    D^p \, B_{L,\, \mathcal{S}}(X) ~\subseteq~
                B_{L,\, \mathcal{S}}(X), \textnormal{ for } p \in \mathbb{N}^n
\end{equation}

In this way we obtain that the foam and multi-foam algebras are \textit{differential algebras
of generalized functions}. \\

Also, the foam and multi-foam algebras contain the Schwartz distributions,
that is, we have the \textit{linear embeddings} which respect the
arbitrary partial derivation of smooth functions

\begin{equation}\label{label1.21}
    \mathcal{D} \, '(X) ~\subset~
             B_{L,\, \Sigma}(X), \textnormal{ for } \Sigma \in \mathcal{S}_{\mathcal{D}}
\end{equation}

\begin{equation}\label{label1.22}
    \mathcal{D} \, '(X) ~\subset~ B_{L,\, \mathcal{S}}(X)
\end{equation}

Indeed, let us recall the wide ranging purely {\it algebraic characterization} of all those
quotient type algebras of generalized functions in which one can embed linearly the Schwartz
distributions, a characterization first given in 1980, see Rosinger [4, pp. 75-88], as well as
Rosinger [5, pp. 306-315], Rosinger [6, pp. 234-244]. According to that characterization -
which also contains the Colombeau algebras as a particular case - the {\it necessary and
sufficient} condition for the existence of the linear embedding (\ref{label1.21}) is precisely
the off diagonality condition in (\ref{label1.13}). Similarly, the necessary and sufficient
condition for the existence of the linear embedding (\ref{label1.22}) is exactly the off
diagonality condition (\ref{label1.15}).

\bigskip
One more property of the foam and multi-foam algebras will prove to be useful. Namely, in view
of (\ref{label1.10}), it is clear that, for every $\Sigma \in \mathcal{S}$, we have the
inclusion $\mathcal{J}_{L,\, \Sigma}(X) \subseteq \mathcal{J}_{L,\, \mathcal{S}}$, and thus we
obtain the \textit{surjective algebra homomorphism}

\begin{equation}\label{label1.23}
    B_{L,\, \Sigma}(X) \ni w + \mathcal{J}_{L,\, \Sigma}(X)
    ~\longmapsto~ w + \mathcal{J}_{L,\, \mathcal{S}}(X) \in B_{L,\, \mathcal{S}}(X)
\end{equation}

And as we shall see in the next subsection, (\ref{label1.23}) can naturally be interpreted as
meaning that the typical generalized functions in $B_{L,\, \mathcal{S}}(X)$ are \textit{more
regular} than those in $B_{L,\, \Sigma}{X}$.

\subsection{Regularity of Generalized Functions}

One natural way to interpret (\ref{label1.23}) in the given context of generalized functions
is the following. Given two spaces of generalized functions $E$ and $F$, such as for instance

\begin{equation}\label{label1.24}
    \mathcal{C}^{\infty}(X) ~\subset~ E ~\subset~ F
\end{equation}

then the larger the space $F$ the \textit{less regular} its typical element can appear to be,
when compared with those of $E$. By the same token, the it smaller the space $E$ , the
\textit{more regular}, compared with those of $F$, one can consider its typical elements. \\

Similarly, given a \textit{surjective} mapping

\begin{equation}\label{label1.25}
    E ~\longrightarrow~ F
\end{equation}

one can again consider that the typical elements of $F$ are at least as \textit{regular} as
those of $E$.

\bigskip
In this way, in view of (\ref{label1.23}), we can consider that, owing to the given
\textit{surjective} algebra homomorphism, the typical elements of the multi-foam algebra
$B_{L,\, \mathcal{S}}(X)$ can be seen as being \textit{more regular} than the typical elements
of the foam algebra $B_{L,\, \Sigma}(X)$.
\\
Furthermore, the algebra $B_{L,\, \mathcal{S}}(X)$ is obtained by factoring the same
$(\mathcal{C}^{\infty}(X))^{\Lambda}$ as in the case of the algebra $B_{L,\, \Sigma}(X)$, this
time however by the significantly \textit{larger} ideal $\mathcal{J}_{L,\, \mathcal{S}_L}(X)$,
an ideal which, unlike any of the individual ideals $\mathcal{J}_{L,\, \Sigma}(X)$, can
simultaneously deal with \textit{all} the singularity sets $\Sigma \in \mathcal{S}_L$, some,
or in fact, many of which can be \textit{dense} in $X$. Further details related to the
connection between \textit{regularization} in the above sense, and on the other hand,
properties of \textit{stability}, \textit{generality} and \textit{exactness} of
generalized functions and solutions can be found in Rosinger [4-6].

\bigskip
This kind of interpretation will be used in section 3 related to the global
Cauchy-Kovalevskaia theorem. Also, it will be further illustrated with the examples of the
differential algebras of generalized functions presented in section 2.

\section{On the Structure of Space-Time Foam Algebras}

\subsection{Special Families of Singularities}

Since in section 3 the space-time foam algebras will be used in order to obtain global
generalized solutions under the usual conditions of the Cauchy-Kovalevskaia theorem, it is
useful to understand the structure of these algebras. And for that, one has to understand the
structure of the ideals $\mathcal{J}_{L,\, \Sigma}(X)$ and $\mathcal{J}_{L,\, \mathcal{S}}(X)$.
In particular, one has to have an idea about their {\it size}. Indeed, in view of the
interpretation at the end of the previous subsection, the \textit{larger} such ideals, the
\textit{more regular} the typical generalized functions in the corresponding quotient algebras
(\ref{label1.11}), (\ref{label1.12}).

\bigskip
In order to be able gain more information about the mentioned algebras, we shall study a more
particular case of them. This case is constructed by allowing a certain relationship between
the right directed partially ordered sets $L=(\Lambda, \leq)$, and the families $\mathcal{S}$
of singularities. Namely, let us take associated with $L$ any set $\mathcal{S}_L$ of subsets
$\Sigma$ of $X$, a set which satisfies the following three conditions. First

\begin{equation}\label{label2.1}
    X \setminus \Sigma \textnormal{ is dense in } X, \,
                     \textnormal{ for } \Sigma \in \mathcal{S}_L
\end{equation}

then, second

\begin{equation}\label{label2.2}
\begin{split}
    &\forall \quad \Sigma, \Sigma \, ' \in \mathcal{S}_L :\\ &
    \exists \quad \Sigma \, '' \in \mathcal{S}_L :\\ &
    \quad \quad \Sigma \cup \Sigma \, ' ~\subseteq~ \Sigma \, ''\\&
\end{split}
\end{equation}

and finally, every $\Sigma \in \mathcal{S}_L$ can be represented as

\begin{equation}\label{label2.3}
    \Sigma ~=~ \limsup_{\lambda \in \Lambda}\, \Sigma_{\lambda} ~=~ \bigcap_{\lambda
    \in \Lambda} \bigcup_{~~~\mu \in \Lambda,~ \mu \geq \lambda}\Sigma_{\mu}
\end{equation}

where $\Sigma_{\lambda}\subseteq X$, while $X\setminus \Sigma_{\lambda}$ is open, for $\lambda
\in \Lambda$.

\bigskip
It is easy to see that we shall have $\mathcal{S}_L \subseteq \mathcal{S}_{\mathcal{D}}$, thus
we are within the framework of the constructions in the previous section.

\bigskip
Further, let us assume that for two subsets $\Sigma,\, \Sigma \, ' \subseteq X$ we have the
representations $\Sigma = \limsup_{\lambda \in \Lambda} \Sigma_{\lambda}$ and $\Sigma\, ' =
\limsup_{\lambda\in \Lambda}\Sigma\, '_{\lambda}$, with $\Sigma_{\lambda},\, \Sigma\, '_{
\lambda} \subseteq X$, where $X \setminus \Sigma_{\lambda}, X \setminus \Sigma\, '_{\lambda}$,
are open, for $\lambda \in \Lambda$. Then, for $\lambda \in \Lambda$, we define $\Sigma\, ''_{
\lambda} = \Sigma_{\lambda} \cup \Sigma \, '_{\lambda}$, hence $X\setminus \Sigma\, ''_{
\lambda}$ is open. In this way, in $X$, the subset $\Sigma\, '' = \limsup_{\lambda \in \Lambda}
\Sigma\, ''_{\lambda}$ has a representation (\ref{label2.3}), and clearly $\Sigma \cup \Sigma
\, ' \subseteq \Sigma\, ''$.
\\
This however, need not mean that (\ref{label2.3}) implies (\ref{label2.2})) since $X \setminus
\Sigma \, ''$ need not be dense in $X$.

\bigskip
We also note that, for a suitable right directed partial order $L = (\Lambda, \leq)$,
condition (\ref{label2.3}) is easy to satisfy for any nonvoid $\Sigma \subseteq X$. Indeed,
let us take as $\Lambda$ the set of all $\lambda = A \subseteq \Sigma$, with nonvoid finite
$A$. Further, for $\lambda = A,\, \mu =B\in \Lambda$, we define the right directed partial
order relation $\lambda \leq \mu$ by the condition $A \subseteq B$. Finally, for $\lambda = A
\in \Lambda$, we take $\Sigma_{\lambda} = A$, in which case relation (\ref{label2.3}) follows
easily.

\bigskip
The above construction shows that in Euclidean spaces, it only has the following three
different cases with respect to the size of $\Lambda$. First, when $\Sigma$ itself is finite.
In this case the above construction can further be simplified, as we can take $\Lambda$ being
composed of one element only, and with the trivial partial order on it, while we take
$\Sigma_{\lambda} = \Sigma$. Then there are the two only other cases, when $\Sigma$ is
countable, respectively, uncountable, and when correspondingly, $\Lambda$ can be taken
countable or uncountable.
\\
Obviously, we may expect to meet in various applications representations (\ref{label2.3})
which are more complicated than those constructed above, at least from the point of view of
the partial orders $\leq$ on $\Lambda$, see for instance Example 1 below.

\bigskip
Let us note that $\mathcal{S}_{nd}$ in (\ref{label1.2}) satisfies (\ref{label2.1}) -
(\ref{label2.3}), if in the last condition we take $L = ( \Lambda, \leq ) = \mathbb{N}$, while
for a given $\Sigma$, and for each $\lambda = \nu \in \Lambda = \mathbb{N}$, we take
$\Sigma_{\lambda} = \Sigma$.

\bigskip
Also, every $\Sigma \in \mathcal{S}_{\textit{Baire I}}$, see (\ref{label1.3}), is of first
Baire category, thus it is a countable union of nowhere dense sets in $X$. In this way
$\mathcal{S}_{\textit{Baire I}}$ satisfies (\ref{label2.1}) - (\ref{label2.3}) for the above
$L$. We also note that the family of singularities $\mathcal{S}_{\textit{Baire I}}$ contains
plenty of singularity sets $\Sigma$ which are \textit{dense} in $X$, and which in addition,
have the cardinal of the continuum.

\subsection{Special Ideals}

Suppose now given a family $\mathcal{S}_L$ of singularities in $X$ satisfying (\ref{label2.1})
- (\ref{label2.3}).
\\
For any singularity set $\Sigma \in \mathcal{S}_L$ and any of its representations in
(\ref{label2.3}), given by a particular family $\mathcal{S}=( \Sigma_{\lambda}\, | \, \lambda
\in \Lambda)$, we denote by

\begin{equation}\label{label2.4}
    \mathcal{I}_{L,\, \Sigma,\, \mathcal{S}}(X)
\end{equation}

the \textit{ideal} in $(\mathcal{C}^{\infty}(X))^{\Lambda}$ consisting of all the sequences of
smooth functions indexed by the set $\Lambda$, namely, $w = ( w_{\lambda} \, | \, \lambda \in
\Lambda ) \in (\mathcal{C}^{\infty}(X))^{\Lambda}$, sequences which \textit{outside} of the
singularity set $\Sigma$ will satisfy the \textit{asymptotic vanishing} condition

\begin{equation}\label{label2.5}
\begin{split}
    & \forall \quad x \in X \setminus \Sigma: \\ &
    \exists \quad \lambda \in \Lambda:\\ &
    \forall \quad \mu \in \Lambda,~ \mu \geq \lambda :\\ &
    \exists \quad x \in \Delta_{\mu} \subseteq X \setminus
                      \Sigma_{\mu},~ \Delta_{\mu} \textnormal{ open } :\\ &
    \quad \quad w_{\mu}=0 \textnormal{ on } \Delta_{\mu}\\ &
\end{split}
\end{equation}

In other words; the sequences of smooth functions $w = ( w_{\lambda} \, | \, \lambda \in
\Lambda )$ in the ideal $\mathcal{I}_{L,\, \Sigma,\, \mathcal{S}}(X)$ are in certain ways
\textit{covering} with their support the singularity set $\Sigma$, while outside of it, they
are vanishing asymptotically.
\\
Also, the ideal $\mathcal{I}_{L,\, \Sigma,\, \mathcal{S}}(X)$ carries in an \textit{algebraic}
manner the information on the singularity set $\Sigma$. It follows that a \textit{quotient} in
which the factorization is made with such ideals may in certain ways \textit{do away with
singularities} through purely algebraic means, see (\ref{label2.11}), (\ref{label2.13}) below.

\bigskip
The ideal $\mathcal{I}_{L,\, \Sigma,\, \mathcal{S}}(X)$ appears to depend not only on $L$ and
$\Sigma$ but also on the family $S = (\Sigma_{\lambda}\, | \, \lambda \in \Lambda)$ which is
in the particular representation of $\Sigma$ in (\ref{label2.3}). However, as we shall see
next in Lemma 1, the ideal $\mathcal{I}_{L,\, \Sigma,\, \mathcal{S}}(X)$ only depends on $L$
and $\Sigma$, and does {\it not} depend on $S$, thus on the representation in (\ref{label2.3}).
Therefore, from now on, this ideal will be denoted by

\begin{equation}\label{label2.6}
    \mathcal{I}_{L,\, \Sigma}(X)
\end{equation}

\textbf{Lemma 1}

\medskip

The ideal $\mathcal{I}_{L,\, \Sigma,\, \mathcal{S}}(X)$ does \textit{not} depend on $S = (
\Sigma_{\lambda}\, | \, \lambda \in \Lambda)$ in the representation (\ref{label2.3}) of
$\Sigma$.

\bigskip
\textbf{Proof.}

\medskip
Let $S = ( \Sigma_{\lambda}\, | \, \lambda \in \Lambda )$ and $S\, ' = ( \Sigma_{\lambda} \, '
\, | \, \lambda \in \Lambda )$ be two representation (\ref{label2.3}) of $\Sigma$. We prove
now that

\begin{equation}\label{label2.7}
    \mathcal{I}_{L,\, \Sigma,\, \mathcal{S}}(X) ~\subseteq~
                     \mathcal{I}_{L,\, \Sigma,\, \mathcal{S}\, '}(X)
\end{equation}

Let us take any sequence of smooth functions $w = ( w_{\lambda} \, | \, \lambda \in \Lambda)$
in the left hand term of (\ref{label2.7}). Then (\ref{label2.5}) holds for $S$. In particular,
for every given $x \in X \setminus \Sigma$, we can find $\lambda \in \Lambda$, such that for
each $\mu \in \Lambda,~ \mu \geq \lambda$, there exists $x \in \Delta_{\mu} \subseteq X
\setminus \Sigma , \, \Delta_{\mu}$ open, and $w_{\mu} = 0$ on $\Delta_{\mu}$.
\\
We now show that by replacing $S$ with $S\, '$, we still have (\ref{label2.5}), for $x \in X
\setminus \Sigma$ arbitrarily given as above. Indeed, in view of (\ref{label2.3})
corresponding to the representation of $\Sigma$ given now by $S\, '$, there exists $\lambda\,'
\in \Lambda$ such that for all $\mu \in \Lambda,~ \mu \geq \lambda \, '$, we have $x \in X
\setminus \Sigma_{\mu} \, '$, and $X \setminus \Sigma_{\mu} \, '$ is by assumption open.
\\
Let us take $\lambda \, '' \in \Lambda$ with $\lambda \, '' \geq \lambda, \, \lambda \, '$.
Then for all $\mu \in \Lambda,~ \mu \geq \lambda \, ''$, we obviously have $x \in \Delta_{\mu}
\, '' \subseteq X \setminus \Sigma_{\mu} \, '$, where $\Delta_{\mu} \, '' = \Delta_{\mu} \cap
(X \setminus \Sigma_{\mu} \, ')$ is open. And clearly, $w_{\mu} = 0$ on $\Delta_{\mu}\, ''$,
since $\Delta_{\mu} \, '' \subseteq \Delta_{\mu}$.
\\
In this way (\ref{label2.5}) does indeed hold for $S\, '$ as well, and the proof of
(\ref{label2.7}) is completed.
\\
Applying (\ref{label2.7}) the other way around, we obtain the proof of Lemma 1.

\qed

\bigskip
Now we establish a few properties of the ideals (\ref{label2.6}). It is easy to see that, for
$\Sigma \in \mathcal{S}_L$, we have

\begin{equation}\label{label2.8}
    \mathcal{I}_{L,\, \Sigma}(X) ~\subseteq~ \mathcal{J}_{L,\, \Sigma}(X)
\end{equation}

also

\begin{equation}\label{label2.9}
    D^p \, \mathcal{I}_{L,\, \Sigma}(X) ~\subseteq~
            \mathcal{I}_{L,\, \Sigma}(X), \textnormal{ for } p \in \mathbb{N}^n
\end{equation}

as well as

\begin{equation}\label{label2.10}
    \mathcal{I}_{L,\, \Sigma}(X) \cap \mathcal{U}_{\Lambda}^{\infty}(X) ~=~ \{~ 0 ~\}
\end{equation}

\subsection{Special Foam Algebras}

Let us consider the \textit{special foam algebras}, Rosinger [9,10], which are defined as
follows for every given $\Sigma \in \mathcal{S}_L$

\begin{equation}\label{label2.11}
    A_{L,\, \Sigma}(X) ~=~ (\mathcal{C}^{\infty}(X))^{\Lambda} / \mathcal{I}_{L,\, \Sigma}(X)
\end{equation}

In view of (\ref{label2.9}), (\ref{label2.10}), these are associative, commutative and unital
differential algebras of generalized functions which contain the Schwartz distributions. And
according to (\ref{label2.8}), for $\Sigma \in \mathcal{S}_L$ we have the \textit{surjective
algebra homomorphism}

\begin{equation}\label{label2.12}
    A_{L,\, \Sigma}(X) \ni w + \mathcal{I}_{L,\, \Sigma}(X) ~\longmapsto~ w +
                                    \mathcal{J}_{L,\, \Sigma}(X) \in B_{L,\, \Sigma}(X)
\end{equation}

And as also in the case of in (\ref{label1.23}), this can be interpreted as indicating that
the typical generalized functions in $B_{L,\, \Sigma}(X)$ are \textit{more regular} than the
typical generalized functions in $A_{L,\, \Sigma}(X)$. More details related to this
interpretation were presented in subsection 1.7.

\subsection{Special Multi-foam Algebras}

Now in order to deal simultaneously with all the singularity sets $\Sigma \in \mathcal{S}_L$,
we define

\begin{equation}\label{label2.13}
    \mathcal{I}_{L,\, \mathcal{S}_L}(X) ~=~ \bigcup_{\Sigma\in \mathcal{S}_L}
                         \mathcal{I}_{L,\, \Sigma}(X)
\end{equation}

which is the \textit{ideal} in $(\mathcal{C}^{\infty}(X))^{\Lambda}$ generated by all ideals
$\mathcal{I}_{L,\, \Sigma} (X)$, with $\Sigma \in \mathcal{S}_L$. Indeed, this is an ideal,
since similar with (\ref{label1.9}), we have the implication in Lemma 2 below.
\\
Finally, we can define the so called \textit{special multi-foam algebra} of generalized
functions, see Rosinger [9,10]

\begin{equation}\label{label2.14}
    A_{L,\, \mathcal{S}_L}(X) ~=~ (\mathcal{C}^{\infty}(X))^{\Lambda} /
                        \mathcal{I}_{L,\, \mathcal{S}_L}(X)
\end{equation}

Together, both the special foam algebras (\ref{label2.11}) and the
special multi-foam algebras in (\ref{label2.14}), will for
simplicity be called \textit{special space-time foam} algebras.

\bigskip
\textbf{Lemma 2.}

\medskip
If~ $\Sigma,~ \Sigma\, ' \in \mathcal{S}_L$, then

\begin{equation}\label{label2.15}
    \Sigma ~\subseteq~ \Sigma\, ' ~~\Longrightarrow~~
                   \mathcal{I}_{L,\, \Sigma}(X) ~\subseteq~ \mathcal{I}_{L,\, \Sigma\, '}(X)
\end{equation}

\textbf{Proof.}

\medskip
Let $\Sigma, \Sigma \, '$ have the corresponding representations in (\ref{label2.3}) given by
$S =(\Sigma_{\lambda}\, | \, \lambda \in \Lambda)$ and $S\, ' =(\Sigma_{\lambda} \, ' \, | \,
\lambda \in \Lambda)$, respectively. Then in view of Lemma 1, the relation (\ref{label2.5})
will now hold for $\Sigma, S,$ and any $w = ( w_{\lambda} \, | \, \lambda \in \Lambda ) \in
\mathcal{I}_{L,\,\Sigma}$. We want to show that (\ref{label2.5}) also holds for $\Sigma \, ',
S \, '$ and $w$. Indeed, let $x \in X\setminus \Sigma \, '$. Then $x \in X \setminus \Sigma$,
since $\Sigma \subseteq \Sigma \, '$. Hence the assumption $w\in \mathcal{I}_{L,\, \Sigma}(X)$,
together with (\ref{label2.5}) give $\lambda \in \Lambda$ such that for every $\mu \in
\Lambda, \, \mu \geq \lambda$ there exists $x \in \Delta_{\mu} \subseteq X \setminus
\Sigma_{\mu}$, with $w_{\mu} = 0$ on $\Delta_{\mu}$ and $\Delta_{\mu}$ is open.
\\
On the other hand, in view of (\ref{label2.3}), we note that $X\setminus \Sigma_{\mu} \, '$ is
open, for $\mu \in \Lambda$. Also $x \in X \setminus \Sigma \, '$ gives $\lambda \, '  \in
\Lambda$ such that $x \in X \setminus \Sigma_{\mu} \, '$, for every $\mu \in \Lambda, \, \mu
\geq \lambda \, '$.  Let us take $\lambda \, '' \in \Lambda, \, \lambda \, '' \geq \lambda,
\lambda \, ' $. It follows that for every $\mu \in \Lambda, \, \mu \geq \lambda \, ''$ we
shall have

$$x \in \Delta_{\mu} \, ' ~=~ \Delta \cap ( X \setminus \Sigma_{\mu} \, ' ) ~\subseteq~
                                        X \setminus \Sigma_{\mu} \, '$$

and clearly, by its above definition, $\Delta_{\mu} \, '$ is open. But $w_{\mu} = 0$ on
$\Delta_{\mu} \, '$, as $\Delta_{\mu} \, ' \subseteq \Delta_{\mu}$.
\\
In this way, indeed, $w\in \mathcal{I}_{L,\, \Sigma \, '}(X)$, and thus the proof of
(\ref{label2.15}) is completed.

\qed

Similar with (\ref{label2.8}) - (\ref{label2.10}), we have the properties

\begin{equation}\label{label2.16}
    \begin{split}
    & \mathcal{I}_{L,\, \Sigma}(X) ~\subseteq~ \mathcal{J}_{L,\, \Sigma}(X)
              ~\subseteq~ \mathcal{J}_{L,\, \mathcal{S}_L}(X), \textnormal{ for } \Sigma \in
    \mathcal{S}_L \\ \\&
    \mathcal{I}_{L,\, \Sigma}(X) ~\subseteq~ \mathcal{I}_{L,\, \mathcal{S}_L}(X)
     ~\subseteq~ \mathcal{J}_{L,\, \mathcal{S}_L}(X), \textnormal{ for } \Sigma \in
    \mathcal{S}_L
    \end{split}
\end{equation}

\begin{equation}\label{label2.17}
    D^p \, \mathcal{I}_{L,\, \mathcal{S}_L}(X) ~\subseteq~
    \mathcal{I}_{L,\, \mathcal{S}_L}(X), \textnormal{ for } p \in \mathbb{N}^n
\end{equation}

\begin{equation}\label{label2.18}
    \mathcal{I}_{L,\, \mathcal{S}_L}(X) \cap \mathcal{U}_{\Lambda}^{\infty}(X) ~=~ \{ 0 \}
\end{equation}

Again, (\ref{label2.17}), (\ref{label2.18}) imply that the special multi-foam algebras
$A_{L,\, \mathcal{S}_L}(X)$ are associative, commutative and unital differential algebras of
generalized functions which contain the Schwartz distributions. Further, in view of
(\ref{label2.16}), for $\Sigma \in \mathcal{S}_L$, we have the \textit{surjective algebra
homomorphism}

\begin{equation}\label{label2.19}
    A_{L,\, \Sigma}(X) \ni w + \mathcal{I}_{L,\, \Sigma}(X) ~~\longmapsto~~
    w + \mathcal{I}_{L,\, \mathcal{S}_L}(X) \in A_{L,\, \mathcal{S}_L}(X)
\end{equation}

as well as the \textit{surjective algebra homomorphism}

\begin{equation}\label{label2.20}
    A_{L,\, \mathcal{S}_L}(X) \ni w + \mathcal{I}_{L,\, \mathcal{S}_L}(X)
    ~~\longmapsto~~ w + \mathcal{J}_{L,\, \mathcal{S}_L}(X) \in B_{L,\, \mathcal{S}_L}(X)
\end{equation}

which together with (\ref{label2.12}), (\ref{label1.23}) will give the commutative diagram of
\textit{surjective} algebra homomorphisms

\begin{equation}\label{label2.21}
\setlength{\unitlength}{1cm}
\thicklines%
\begin{picture}(5,5)(2.8,2.5)
\put(0,5){$A_{L,\, \Sigma}(X)$}%
\put(1.7,5.4){\vector(2,1){1.8}}%
\put(1.7,4.8){\vector(2,-1){1.8}}%
\put(3.8,6.5){$B_{L,\, \Sigma}(X)$}%
\put(3.8,3.6){$A_{L,\, \mathcal{S}_L}(X)$}%
\put(5.7,6.3){\vector(2,-1){1.8}}%
\put(5.7,3.8){\vector(2,1){1.8}}%
\put(7.8,5){$B_{L,\, \mathcal{S}_L}(X)$}%
\end{picture}
\end{equation}

In view of subsection 1.7., the interpretation of (\ref{label2.21}) is that the typical
generalized function in the algebras which are the target of arrows are \textit{more regular}
than those in the algebras which are the source of the arrows. In particular the \textit{most
regular} differential algebras of generalized functions among those constructed in this paper
are the \textit{multi-foam algebras} $B_{L,\, \mathcal{S}_L}(X)$.

\subsection{Special Space-Time Foam Algebras}

Again, it will on occasion be convenient to call the special foam algebras and the special
multi-foam algebras by one single term, namely, special space-time foam algebras.

\bigskip
As an illustration of the above, let us recall the \textit{nowhere dense} differential
algebras of generalized functions $A_{nd}(X)$ introduced in Rosinger [3], see also Rosinger
[4-8], and section 3 below, which were recently used in Mallios \& Rosinger [1] as the
structure coefficients replacing the smooth functions in the abstract differential geometry
developed in Mallios [1,2]. \\
Namely,if we take as our right directed partial order the natural numbers $\mathbb{N}$, that
is, $L = ( \Lambda, \leq) = \mathbb{N}$, while we take the family of singularities in $X$ as
given by $\mathcal{S}_{nd}$ in (\ref{label1.2}), then it is easy to see that

\begin{equation}\label{label2.22}
    A_{nd}(X) ~=~ B_{\mathbb{N},\, \mathcal{S}_{nd}}(X) ~=~
                                 A_{\mathbb{N},\, \mathcal{S}_{nd}}(X)
\end{equation}

This, in particular, further clarifies the extent to which Mallios \& Rosinger [2,3]
strengthens the results in Mallios \& Rosinger [1]. Indeed in the latter paper the singularity
sets $\Sigma \subset X$ had to be closed and nowhere dense in $X$, thus their complementaries
were open and dense in $X$. On the other hand, in Mallios \& Rosinger [2,3], we only ask that
for each singularity set $\Sigma \subset X$, the corresponding set of nonsingular points $X
\setminus \Sigma$ be dense in $X$. \\
In this regard, the best previous results in the literature only allowed one single
singularity set $\Sigma \subset X$, which in addition had to be closed and nowhere dense and
in the boundary of $X$, see Heller [2], Heller \& Sasin [2]. Consequently, already the result
in Mallios \& Rosinger [1] - which is considerably more particular than that in Mallios \&
Rosinger [2,3] - proved to be significantly more powerful, since it could deal simultaneously
with {\it all} closed and nowhere dense singularity sets $\Sigma \subset X$, be they in the
boundary of $X$ or not.

\subsection{On the Structure of the Ideals}

Let us show in Example l next, that for \textit{large} enough - that is, uncountable and of
the cardinality of the continuum - index sets $\Lambda$, the ideals $\mathcal{I}_{L,\, \Sigma}
(X)$ are \textit{not} trivial, namely, the do {\it not} collapse to the null ideal $\{~ 0 ~\}$.
For that purpose, for any singularity set $\Sigma \in {\cal S}_{\cal D}$, and for suitably
chosen index sets $\Lambda$, we shall construct \textit{large classes} of sequences of smooth
functions

$$ w^* ~=~ (w^*_{\lambda} \, | \, \lambda \in \Lambda) \in
                                     (\mathcal{C}^{\infty}(X))^{\Lambda} $$

such that $w^* \in \mathcal{I}_{L,\, \Sigma} (X)$. Then in view of (\ref{label2.16}), this
will suffice to show that none of the other ideals $\mathcal{J}_{L,\, \Sigma}(X),
\mathcal{J}_{L,\, \mathcal{S}}(X), \mathcal{I}_{L,\, \mathcal{S}_L}(X)$ is trivial. \\

The construction in Example 1 next will also illustrate in more detail the way \textit{dense}
singularities can be dealt with by purely {\it algebraic} means. \\

Here we should note that the issue of nontriviality of these ideals is in itself nontrivial.
Indeed, in view of the off diagonality conditions (\ref{label1.15}), or equivalently, of the
algebra embeddings (\ref{label1.16}), none of the ideals $\mathcal{J}_{L,\, \Sigma}(X),
\mathcal{J}_{L,\, \mathcal{S}}(X), \\ \mathcal{I}_{L,\, \Sigma}(X)$ or $\mathcal{I}_{L,\,
\mathcal{S}_L}(X)$ can be too large. Thus the issue of the nontriviality of these ideals
involves a conflict. \\
As for determining which are the {\it maximal} ones among such ideals this is still an open
problem, and one that has an obvious importance, as argued in some detailed in Rosinger [15].

\bigskip
\textbf{Example 1.}

\medskip
Given any nonvoid singularity $\Sigma \subset X$ such that $\Sigma \in {\cal S}_{\cal D}$, in
other words, for which $X\setminus \Sigma$ is dense in $X$, let us take the index set
$\Lambda$ as the set of elements

$$ \lambda ~=~ ( A, (\alpha_x \, | \, x\in \Sigma ) ) $$

where $A \subseteq \Sigma$ is nonvoid finite, and for $x \in \Sigma$, we have $\alpha_x \in
\mathcal{D}(X), \, \alpha_x \neq 0$. Here we recall that $\mathcal{D}(X)$ denotes the space of
compactly supported smooth functions on $X$. \\
Now we define on $\Lambda$ a {\it right directed} partial order $\leq$ as follows. Given
$\lambda = (A, (\alpha_x \, | \, x \in \Sigma)), \, \mu = (B, (\beta_x \, | \, x \in \Sigma))$,
we shall write $\lambda \leq \mu$, if and only if

\begin{equation}\label{label2.23}
    A \subseteq B~~~ \textnormal{ and }~~~
           \bigcup_{x \in \Sigma} \textnormal{supp } \beta_x ~\subseteq~
           \bigcup_{x\in \Sigma} \textnormal{supp } \alpha_x
\end{equation}

Further, aiming to obtain for $\Sigma$ a representation (2.3) , for $\lambda = (A, (\alpha_x
\, | \, x \in \Sigma)) \in \Lambda$, we define $\Sigma_{\lambda} = A$. Finally, we also define
the compactly supported smooth function

\begin{equation}\label{label2.24}
    w_{\lambda}^* ~=~ \sum_{x\in A} \alpha_x \in \mathcal{D}(X)
\end{equation}

which is well defined since $A$ is a nonvoid finite set.
\\
Then it follows easily that $X\setminus \Sigma_{\lambda}$ is open for $\lambda \in \Lambda$,
and in addition we also have

\begin{equation}\label{label2.25}
    \Sigma ~=~ \limsup_{\lambda \in \Lambda} \Sigma_{\lambda}
\end{equation}

The fact of interest to us is that

\begin{equation}\label{label2.26}
    w^* ~=~ (w^*_{\lambda}\, | \, \lambda \in \Lambda) \in \mathcal{I}_{L,\, \Sigma}(X)
\end{equation}

Indeed, for the proof of (\ref{label2.26}), let us take any $y\in
X\setminus \Sigma$, then (\ref{label2.25}) gives

\begin{align*}
\exists \quad \lambda &= (A, (\alpha_x \, | \, x\in \Sigma)) \in \Lambda  :\\
\forall \quad \mu & = (B, (\beta_x \, | \, x\in \Sigma)) \in \Lambda,~ \mu \geq \lambda  :\\
\quad \quad y & \in X \setminus \Sigma_{\mu} ~=~ X \setminus B
\end{align*}

But obviously, we can assume that

\begin{equation}\label{label2.27}
    y \not\in \bigcup_{x \in \Sigma}~ \textnormal{supp}\, \,  \alpha_x
\end{equation}

since we took $y \in X \setminus \Sigma$, while for $x\in \Sigma$, we have $\alpha_x \in
\mathcal{D}(X)$, and the only other condition $\alpha_x$ has to satisfy is that $\alpha_x(x)
\neq 0$.

\medskip
Now let $\mu = (B, (\beta_x\, | \, x \in \Sigma)) \in \Lambda,~ \mu \geq \lambda$, then $A
\subseteq B$ and

\begin{equation*}
\bigcup_{x \in \Sigma} \textnormal{supp} \, \, \beta_x ~\subseteq~
\bigcup_{x \in \Sigma} \textnormal{supp} \, \, \alpha_x
\end{equation*}

hence the previous assumption (\ref{label2.27}) gives

$$y\not\in \bigcup_{x \in \Sigma} \textnormal{supp } \beta _x$$

thus

$$\exists \, \Delta_{\mu} \textnormal{ open, }~ y \in \Delta_{\mu}
                \subseteq X \setminus \Sigma_{\mu} = X \setminus B \, :$$

$$\Delta_{\mu}\bigcap ~(~ \bigcup_{x\in B}\textnormal{ supp }\beta_x ~) ~=~ \phi$$

In this way we obtain that $w^*_{\mu} = \sum_{x\in B}\beta_x=0$ on $\Delta_{\mu}$, and in view
of (\ref{label2.5}), the proof of (\ref{label2.25}) is completed.

\qed

\bigskip
From the point of view of dealing with \textit{dense} singularities, the essential property in
Example 1 above is illustrated in the sequences of smooth functions constructed in
(\ref{label2.26}), namely

\begin{equation}\label{label2.28}
    w^* ~=~ (w_{\lambda}^^ \, |\, \lambda \in \Lambda)\in \mathcal{I}_{L,\, \Sigma}(X)
\end{equation}

Indeed, in view of (\ref{label2.27}), these sequence have the property

\begin{equation}\label{label2.29}
    \begin{split}
    & \forall \quad y \in X\setminus\Sigma : \\ &
    \exists \quad \alpha_x \in \mathcal{D}(X), \textnormal{ for each } x\in \Sigma :\\ &
    \quad \quad y \not\in \bigcup_{x\in \Sigma} \textnormal{ supp } \alpha_{x}
    \end{split}
\end{equation}

which means that

\begin{equation}\label{label2.30}
    \begin{split}
    &\forall \quad x \in X \setminus \Sigma \, :\\ &
    \exists \quad \lambda \in \Lambda \, :\\&
    \forall \quad \mu \in \Lambda,~ \mu \geq \lambda \, :\\&
    \quad \quad x \not \in \textnormal{supp } w^*_{\mu}
    \end{split}
\end{equation}

On the other hand, owing to the specific definition of $\Lambda$ in Example 1, it follows that

\begin{equation}\label{label2.31}
    \begin{split}
    & \forall \quad \lambda = (A, (\alpha_x\, | \, x\in \Sigma))\in \Lambda \, :\\&
    \quad \quad \phi \neq \Sigma ~\subseteq~ \bigcup_{x\in \Sigma}
                                          \textnormal{ supp } \alpha_x
    \end{split}
\end{equation}

Here, as a consequence, we note \textit{four} facts related to the sequences in
(\ref{label2.28}).

\begin{itemize}

\item First, in (\ref{label2.30}) the \textit{singularity} sets $\Sigma \subset X$ can be
arbitrary large, provided that their complementary $X \setminus \Sigma$ are still dense in $X$.
In particular, $\Sigma$ can have the cardinal of the continuum while $X \setminus \Sigma$ need
only be countable and dense in $X$. As we mentioned, in the case of the real line $X =
\mathbb{R}$, for instance, $\Sigma$ can be the uncountable set of all irrational numbers,
since its complementary $X\setminus \Sigma$, that is, the rational numbers, is still dense in
$X$, although it is only countable.
\\
And yet, every point $x$ outside of such rather arbitrary singularity sets $\Sigma$ will
eventually also be outside of the support of $w_{\lambda}^*$, see earlier comment following
(\ref{label1.8}).

\item Second, due to (\ref{label2.31}), such rather arbitrary singularity sets $\Sigma \subset
X$ will nevertheless be included in the support of the functions $\alpha_x$ which through
(\ref{label2.24}), make up the terms of the sequences in (\ref{label2.28}), sequences which
guarantee the nontriviality of the mentioned ideals.

\item Third, the \textit{index} set $\Lambda$ can \textit{depend} on the given singularity set
$\Sigma$, and can be rather \textit{large}. In particular, it may even have to be uncountable
and of the cardinality of the continuum, as happens in Example 1 above, and for the
corresponding sequences in (\ref{label2.28}).

\item Finally, we note that the above can give us a certain information about the possible
\textit{size} of the various deals we have considered so far.  Indeed, in view of
(\ref{label2.16}), we obtain for every singularity set $\Sigma \in \mathcal{S}_L$

\begin{equation}\label{label2.32}
    \begin{split}
    & w^* \in \mathcal{I}_{L,\, \Sigma}(X) ~\subseteq~ \mathcal{J}_{L,\, \Sigma}(X)
                 ~\subseteq~ \mathcal{J}_{L,\, \mathcal{S}_L}(X)\\ \\&
    w^* \in \mathcal{I}_{L,\, \Sigma}(X) ~\subseteq~ \mathcal{I}_{L,\, \mathcal{S}_L}(X)
                ~\subseteq~ \mathcal{J}_{L,\, \mathcal{S}_L}(X)
    \end{split}
\end{equation}

\end{itemize}

\section{Global Cauchy-Kovalevskaia Theorem}

\subsection{Preliminary Comments}

Let us recall that, as mentioned, in Rosinger [1-8] a large variety, and in fact,
\textit{infinitely} many classes of differential algebras of generalized functions were
constructed.
\\
And in case these constructions start with Banach algebra valued, and not merely with real or
complex valued functions, the resulting algebras can be \textit{noncommutative} as well. Also
a wide ranging and purely algebraic characterization was given in Rosinger [4-6] for those
algebras which contain the linear vector space of Schwartz distributions. That
characterization is expressed by the {\it off diagonality} condition whose specific instance
can be seen in (1.15), for instance. And in view of that characterization, the Colombeau
algebras prove to be a particular case of the algebras constructed earlier in Rosinger [1-4],
see details in Rosinger [5,6] or Grosser et.al. [p. 7]. \\

Until recently, only two particular cases of these classes of algebras have been used in the
study of global generalized solutions of nonlinear PDEs. Namely, first was the class of the
nowhere dense differential algebras of generalized functions, see Rosinger [1-8], or
(\ref{label2.22}) above, while later came the class of algebras considered in Colombeau.
\\
These latter algebras, since they also contain the Schwartz distributions are - in view of the
above mentioned algebraic characterization - by necessity a particular case of the classes of
algebras of generalized functions first introduced in Rosinger [1-8]. \\

The Colombeau algebras of generalized functions enjoy a rather simple and direct connection
with the Schwartz distributions, and therefore, with a variety of Sobolev spaces as well.
Furthermore, the polynomial type growth conditions which define - and also seriously limit
with respect to the singularities which they are able to handle - the generalized functions in
the Colombeau algebras, can offer an easy and familiar set up to work with for certain
analysts. This led to their relative popularity in the study of generalized solutions of
PDEs. \\
What happens, however, is that the ease and familiarity of working with growth type conditions
not only restricts the class of singularities which can be dealt with, but also leads quite
soon to considerable technical complications. \\
One such instance can be seen when comparing the difficulties in defining on finite
dimensional smooth manifolds the Colombeau algebras, see for instance, Grosser et.al., and on
the other hand, the rather immediate and natural manner in which space-time foam algebras can
be defined on the same manifolds, see Rosinger [11]. \\

Another instance, related to the main subject of this paper, and already mentioned, is the
following. In view of the severe limitations on the class of singularities the Colombeau
algebras of generalized functions are able to deal with, one simply cannot formulate, let
alone obtain in such algebras the global version of the Cauchy-Kovalevskaia theorem presented
in this paper. \\

However, as mentioned in Rosinger [7, pp. 5-8, 11-12, 173-187], one should avoid rushing into
a too early normative judgement about the way the long established \textit{linear} theories of
generalized functions - such as for instance the Schwartz or Sobolev distributions - should
relate to the still emerging, and far more complex and rich corresponding {\it nonlinear}
theories. \\
In particular, \textit{two} aspects of such possible relationships still await a more
thoroughly motivated and clear settlement :

\begin{itemize}

\item First, the purely algebraic-differential type connections between Schwartz distributions
and the more recently constructed variety of differential algebras of generalized functions
should be studied in more detail. And since the main aim is to deal with generalized - hence,
not smooth enough, but rather singular - solutions, a main stress should be placed on the
respective capabilities to deal with singularities, see some related comments in Rosinger [8,
pp. 174-185]. In this regard, let us only mention the following. \\
A fundamental property of various spaces of generalized functions which is closely related to
their capability to handle a large variety of singularities is that such spaces should have a
\textit{flabby sheaf} structure, see for details Kaneko, or Rosinger [10,13]. However, the
Schwartz or Sobolev distributions, the Colombeau generalized functions, as well as scores of
other frequently used spaces of generalized functions happen to {\it fail} being flabby
sheaves. \\
On the other hand, the nowhere dense differential algebras of generalized functions, see
(\ref{label2.22}) above, have a structure of flabby sheaves, as shown for instance in Mallios
\& Rosinger [1]. \\
Similarly, the far larger class of space-time foam differential algebras of generalized
functions dealt with in this paper prove to be flabby sheaves as well, see Rosinger [11] or
Mallios \& Rosinger [2,3].

\item The second aspect is possibly even more controversial. And it is so, especially because
of the historical phenomenon that the study of the \textit{linear} theories of generalized
functions has from its early modern stages in the 1930s been strongly connected with the then
massively emerging theories of linear topological structures.
\\
However, just as with the \textit{nonstandard} reals $^*\mathbb{R}$, so with the various
differential algebras of generalized functions, it appears that \textit{infinitesimal} type
elements in these algebras play an important role. And the effect is that if one introduces
Hausdorff topologies on these algebras, then, when these topologies are restricted to the more
regular, smooth, classical type functions, they inevitably lead to the trivial
\textit{discrete} topology on them, see related comments in the mentioned places in Rosinger
[8], as well as in Remark 1 below.

\end{itemize}

Compared, however, with the nowhere dense differential algebras of generalized functions, let
alone the space-time foam ones dealt with in this paper, the Colombeau algebras suffer from
several important limitations. Among them, relevant to this paper is the following.

\begin{itemize}

\item There are polynomial type \textit{growth conditions} which the generalized functions must
satisfy in the neighbourhood of singularities.

\end{itemize}

On the other hand, the earlier introduced nowhere dense algebras do not suffer from any of the
above two limitations. Indeed, the nowhere dense algebras allow singularities on arbitrary
closed nowhere dense sets, therefore, such singularity sets can have arbitrary large positive
Lebesgue measure, Oxtoby. Furthermore, in the nowhere dense algebras no any kind of conditions
are asked on generalized functions in the neighbourhood of their singularities.
\\
In fact, it is precisely due to the lack of the mentioned type of constraints that the nowhere
dense algebras have a flabby sheaf structure, while the Colombeau algebras fail to do so. \\

Here, for the sake of clarity, let us briefly elaborate on the above. The space ${\cal
C}^\infty ( \mathbb{R}^n )$ is of course a subset of the Colombeau algebra on $\mathbb{R}^n$.
However, the smallest flabby sheaf containing ${\cal C}^\infty ( \mathbb{R}^n )$ is, [34, pp.
143-146] \\

$~~~~~~ {\cal C}^\infty_{nd} ( \mathbb{R}^n ) = \left \{~ f : \mathbb{R}^n \longrightarrow
               \mathbb{C} ~~ \begin{array}{|l} ~~
                  \exists~~ \Gamma \subset \mathbb{R}^n,~
                    \mbox{closed, nowhere dense} ~: \\ \\
                      ~~~~~~ f|_{\mathbb{R}^n \setminus \Gamma}
                     \in {\cal C}^\infty ( \mathbb{R}^n \setminus \Gamma )
                                 \end{array} ~ \right \} $ \\ \\

and this set of functions is no longer contained in the Colombeau algebra on $\mathbb{R}^n$,
since that algebra fails to be a flabby sheaf. More precisely, an arbitrary function
$f|_{\mathbb{R}^n \setminus \Gamma} \in {\cal C}^\infty ( \mathbb{R}^n \setminus \Gamma )$,
where $\Gamma$, for instance, has a positive Lebesgue measure, can have singularities which
the Colombeau algebra on $\mathbb{R}^n$ cannot deal with, since on $\mathbb{R}^n \setminus
\Gamma$ and in the neighbourhood of $\Gamma$, such a function $f|_{\mathbb{R}^n \setminus
\Gamma} \in {\cal C}^\infty ( \mathbb{R}^n \setminus \Gamma )$ can grow far faster than any
polynomial. \\

In this paper, the use of the \textit{space-time foam} and differential algebras of
generalized functions, introduced recently in Rosinger [9,10], brings a further significant
enlargement of the possibilities already given by the nowhere dense algebras. Indeed, this
time the singularities can be concentrated on \textit{arbitrary} subsets, including
\textit{dense} ones, provided that their complementary, that is, the set of nonsingular points,
is still dense. Furthermore, as already in the case of the nowhere dense algebras, also in the
space-time foam algebras, \textit{no} any sort of condition is asked on the generalized
functions in the neighbourhood of their singularities. \\

One interest in obtaining solutions which are in the space-time foam algebras is that, in view
of the interpretations in subsection 1.7., such solutions can be seen as having better
\textit{regularity} properties than those obtained earlier in the nowhere dense algebras. \\

We recall that one could already obtain in the framework of the nowhere dense algebras a very
general, and in fact, \textit{type independent} version of the classical Cauchy-Kovalevskaia
theorem, see Rosinger [4-8], according to which every analytic nonlinear PDE, together with
every non-characteristic analytic initial value problem has a \textit{global} generalized
solution, which is \textit{analytic} on the whole domain of definition of the respective PDE,
except for a closed nowhere dense set, set which if so desired, can be chosen to have zero
Lebesgue measure.
\\
This earlier global and type independent existence results is, fortunately, preserved in the
case of the multi-foam algebras as well, and as seen next, it is significantly {\it
strengthened} with respect to the {\it regularity} properties of the respective global
solutions. \\

Here it can be mentioned that, so far, one could not obtain any kind of similarly general,
powerful, and in fact, type independent existence of global solutions result in any of the
infinitely many other classes of algebras of generalized functions, including in the
Colombeau class of algebras. \\

And as also mentioned, as far as the Colombeau algebras are concerned, they do not allow even
the mere formulation of the global Cauchy-Kovalevskaia theorem, let alone its solution, as
obtained in this paper. \\

\subsection{The Global Cauchy-Kovalevskaia Theorem}

We shall present now in the framework of multi-foam differential algebras of generalized
functions the corresponding global version of the Cauchy-Kovalevskaia theorem. \\

First however, for convenience, let us recall this classical local theorem in its first global
formulation, as it was given for the nowhere dense differential algebras of generalized
functions, see Rosinger [4-10]. \\

We consider the \textit{general nonlinear analytic} partial differential operator

\begin{equation}\label{label3.1}
    T(x, D)U(x) ~=~ D^m_t U(t,y) - G(t,y,...,D^p_t D^q_y U(t,y),...)
\end{equation}

where $U : X \longrightarrow \mathbb{C}$ is the unknown function, while $x = (t,y) \in X,~
t \in \mathbb{R},~ y \in \mathbb{R}^{n-1},~ p \in \mathbb{N},~ 0 \leq p < m,~ q \in
\mathbb{N}^{n-1},~ p+|q| \leq m$, and $G$ is arbitrary analytic in all of its variables. \\

Now together with the analytic nonlinear PDE

\begin{equation}\label{label3.2}
    T(x, D)U(x) ~=~ 0,~~~ x \in X
\end{equation}

we consider the non-characteristic analytic hypersurface

\begin{equation}\label{label3.3}
    S ~=~ \{~ x ~=~ (t,y)\in X ~~|~~ t ~=~ t_0 ~\}
\end{equation}

for any given $t_0\in \mathbb{R}$, and on it, we consider the initial value problem

\begin{equation}\label{label3.4}
    D^p_t U (t_0, y) ~=~ g_p (y),~~~ 0 \leq p < m,~~ (t_0,y)\in S
\end{equation}

Obviously, the analytic nonlinear partial differential operator $T(x, D)$ in (\ref{label3.1})
generates a mapping

\begin{equation}\label{label3.5}
    T(x, D) : \mathcal{C}^{\infty}(X) ~\longrightarrow~ \mathcal{C}^{\infty}(X)
\end{equation}

also, in view of (\ref{label1.18}), (\ref{label1.19}), (\ref{label2.22}) and (\ref{label3.1}),
it generates a mapping

\begin{equation}\label{label3.6}
    T(x, D): A_{nd}(X) ~\longrightarrow~ A_{nd}(X)
\end{equation}

and the mappings (\ref{label3.5}), (\ref{label3.6}), (\ref{label1.16}) form a commutative
diagram

\begin{equation}\label{label3.7}
\setlength{\unitlength}{1cm}
\thicklines%
\begin{picture}(5,3.5)(1.0,1.5)
\put(0,4){$\mathcal{C}^{\infty}(X)$}%
\put(1.3,4.1){\vector(1,0){3}}%
\put(2.1,4.3){$T(x,D)$}%
\put(4.6,4){$\mathcal{C}^{\infty}(X)$}%
\put(5,3.8){\vector(0,-1){1.5}}%
\put(1.4,2){\vector(1,0){2.9}}%
\put(2.1,1.5){$T(x,D)$}%
\put(4.6,1.9){$A_{nd}(X)$}%
\put(0,1.9){$A_{nd}(X)$}%
\put(0.4,3.8){\vector(0,-1){1.5}}%
\end{picture}
\end{equation}

\bigskip
In this way, see Rosinger [4-8], we could obtain the earlier global existence result in the
nowhere dense differential algebras of generalized functions \\ \\

\textbf{Theorem G C-K}

\medskip
The analytic nonlinear PDE in (\ref{label3.2}), with the analytic non-characteristic initial
value problem (\ref{label3.3}), (\ref{label3.4}) has \textit{global} generalized solutions

\begin{equation}\label{label3.8}
    U \in A_{nd}(X)
\end{equation}

defined on the whole of $X$. These solutions $U$ are \textit{analytic} functions

\begin{equation}\label{label3.9}
    \psi : X\setminus \Sigma ~\longrightarrow~ \mathbb{C}
\end{equation}

when restricted to the \textit{open dense} subsets $X \setminus \Sigma$, where the singularity
subsets

\begin{equation}\label{label3.10}
    \Sigma ~\subset~ X,~~~ \Sigma ~\textnormal{ closed, nowhere dense in } X
\end{equation}

can be suitably chosen. Further, one can choose $\Sigma$ to have zero Lebesgue measure, namely

\begin{equation}\label{label3.11}
    \textnormal{mes }\Sigma ~=~ 0
\end{equation}

\qed

\bigskip
As a main result of this paper, we shall \textit{strengthen} the above global existence
theorem by showing that it also holds in certain classes of multi-foam differential algebras
of generalized functions. This will indeed be a strengthening since, as shown next, the
respective classes of algebras $A_{\textit{Baire I}}(X)$ and $B_{\textit{Baire I}}(X)$ in
(3.12), (3.13) below

\begin{itemize}

\item are \textit{surjective} images through algebra homomorphisms of the nowhere dense
algebras used in Theorem G C-K, and furthermore, they

\item are significantly {\it smaller} than the nowhere dense algebras.

\end{itemize}

In this way, since now we shall obtain the existence of global solutions in \textit{smaller}
algebras, this result can also be interpreted as a \textit{regularity} result which improves
on the earlier existence result in Theorem G C-K, see the respective interpretations in
subsection 1.7. \\

Let us for that purpose return to the two classes of singularities in (\ref{label1.2}) and
(\ref{label1.3}), namely $\mathcal{S}_{nd}$ and $\mathcal{S}_{\textit{Baire I}}$, respectively.
Further, as in (\ref{label2.22}), let us take for both of them the same right directed partial
order, given by $L = ( \Lambda, \leq) = \mathbb{N}$. \\

For simplicity, let us denote by $\mathcal{J}_{nd}(X)$ and $\mathcal{J}_{\textit{Baire I}}(X)$
the respective ideals (\ref{label1.10}) which correspond to these two classes of singularities.
Similarly, let us denote by $\mathcal{I}_{nd}(X)$ and $\mathcal{I}_{\textit{Baire I}}(X)$ the
respective ideals (\ref{label2.13}) which correspond to the same two classes of
singularities. \\

Then (\ref{label2.22}) gives the nowhere dense algebra $A_{nd}(X)$ both as the multi-foam
algebra (\ref{label1.12}) defined by the ideal $\mathcal{J}_{nd}(X)$, as well as the special
multi-foam algebra (\ref{label2.14}) defined by the ideal $\mathcal{I}_{nd}(X)$. \\

Let us now denote by $A_{\textit{Baire I}}(X)$ the special multi-foam algebra which in a
similar way is defined by the ideal $\mathcal{I}_{\textit{Baire I}}(X)$. \\

Further, (\ref{label1.12}) similarly gives the multi-foam algebra $B_{\textit{Baire I}}(X)$
defined by the ideal $\mathcal{J}_{\textit{Baire I}}(X)$. \\

Now in view of (\ref{label1.4}), (\ref{label1.10}), (\ref{label2.8}) and (\ref{label2.13}), it
is clear that

\begin{equation}\label{label3.12}
    \mathcal{I}_{nd}(X) ~\subset~ \mathcal{I}_{\textit{Baire I}}(X) ~\subset~
              \mathcal{J}_{\textit{Baire I}}(X)
\end{equation}

This obviously means that we have the \textit{surjective} algebra homomorphisms

\begin{equation}\label{label3.13}
    \begin{split}
    & ~~~~~~A_{nd}(X) ~~~~~~~~\longrightarrow~~~~~~~ A_{\textit{Baire I}}(X)
                   ~~~~~~\longrightarrow~~~~~~~~ B_{\textit{Baire I}}(X)\\&
    U ~=~ s + \mathcal{I}_{nd}(X) ~\longmapsto~ U_{*} =
            s  + \mathcal{I}_{\textit{Baire I}}(X) ~\longmapsto~ U_{**} ~=~
                        s + \mathcal{J}_{\textit{Baire I}}(X)
    \end{split}
\end{equation}

which commute with arbitrary partial derivatives, see (\ref{label1.20}), (\ref{label2.17}).
And in view of the interpretations in subsection 1.7., we should recall that (\ref{label3.13})
means that the typical generalized functions in $B_{\textit{Baire I}}(X)$ are \textit{more
regular} than those both in $A_{nd}(X)$ and $A_{\textit{Baire I}}(X)$. \\

In this way, as a main result of this paper, we obtain the following \textit{global}
Cauchy-Kovalevskaia existence result in algebras with {\it dense singularities}, a result
which also gives \textit{better regularity} properties than those known so far, namely, in the
above mentioned earlier Theorem G C-K, see Rosinger [7] : \\

\textbf{Theorem 1.}

\medskip
The result in Theorem G C-K above holds in any of the following two stronger forms as far as
the {\it regularity} of global solutions is concerned, namely, with

\begin{equation}\label{label3.14}
    U \in A_{\textit{Baire I}}(X)
\end{equation}

or with

\begin{equation}\label{label3.15}
    U \in B_{\textit{Baire I}}(X)
\end{equation}

\textbf{Proof.}

\medskip

In the proof of Theorem G C-K, the global generalized solution $U \in A_{nd}(X)$ is obtained
as given by, see Rosinger [4-7]

$$ U ~=~ s + \mathcal{I}_{nd}(X) \in A_{nd}(X) $$

where

$$ s ~=~ (\psi_{\nu}\, | \, \nu \in \mathbb{N}) \in (\mathcal{C}^{\infty}(X))^{\mathbb{N}} $$

and for every compact $K \subset X \setminus \Sigma$ there exists $\nu \in \mathbb{N}$, such
that $\psi_{\nu}=\psi$ on $K$, for $\mu \in \mathbb{N}, \, \mu \geq \nu$.
\\
Now according to the surjective algebra homomorphisms (\ref{label3.13}), we can take the same
sequence $s$ and define

\begin{equation}\label{label3.16}
    U = s + \mathcal{I}_{\textit{Baire I}} \in A_{\textit{Baire I}}(X)
\end{equation}

or alternatively

\begin{equation}\label{label3.17}
    U = s+ \mathcal{J}_{\textit{Baire I}}\in B_{\textit{Baire I}}(X)
\end{equation}

Then clearly, (\ref{label3.16}), (\ref{label3.17}) will give
respectively the global solutions in (\ref{label3.14}) and
(\ref{label3.15}).

\subsection{ Connections with Distributions}

Let us indicate in short the way the multi-foam algebras can be related to the Schwartz
distributions.  For that, let us recall the wide ranging and purely algebraic characterization,
mentioned at the beginning of this section, of all those differential algebras of generalized
functions in which one can embed linearly the Schwartz distributions, a characterization which,
as also mentioned, contains the Colombeau algebras as a particular case, see Rosinger [4, pp.
75-88], Rosinger [5, pp. 306-315], Rosinger [6, pp. 234-244]. \\

According to the mentioned characterization, in the case of the multi-foam algebras $B_{L,\,
\mathcal{S}}(X)$, for instance, the necessary and sufficient condition for the existence of
such a linear embedding

\begin{equation}\label{label3.18}
    \mathcal{D}\, '(X) ~\subset~ B_{L,\, \mathcal{S}}(X)
\end{equation}

is precisely the off diagonality condition (\ref{label1.15}), which as we have seen, does
indeed hold. And the linear embedding (\ref{label3.18}) will preserve the differential
structure of $\mathcal{C}^{\infty}(X)$. \\

In a similar way, in the case of the special multi-foam algebras, the corresponding off
diagonality condition (\ref{label2.18}) is again the necessary and sufficient condition for
the existence of the linear embedding

\begin{equation}\label{label3.19}
    \mathcal{D}\, '(X) \subset A_{L,\, \mathcal{S}_L}(X)
\end{equation}

which again will preserve the differential structure of $\mathcal{C}^{\infty}(X)$.

\subsection{Final Remarks}

\textbf{Remark 1.}

\medskip
It is important to note that, just like in Mallios \& Rosinger [1], where the nowhere dense
differential algebras of generalized functions were used, or for that matter in Rosinger
[1-11], or Colombeau, Biagioni, Oberguggenberger, Grosser et.al., where other differential
algebras of generalized functions appeared as well, so in this paper, where the space-time
foam and special space-time foam differential algebras of generalized functions are employed,
there is again \textit{no need} for any topological algebra structure on these algebras.
\\
One of the reasons for the lack of need for any topological algebra structure on the algebras
of generalized functions under consideration is the following. It is becoming more and more
clear that the classical Kuratowski-Bourbaki topological concept is not suited to the
mentioned algebras of generalized functions. Indeed, these algebras prove to contain {\it
nonstandard} type of elements, that is, elements which in a certain sense are infinitely small,
or on the contrary, infinitely large. And in such a case, just like in the much simpler case
of nonstandard reals $^*\mathbb{R}$, any topology which would be Hausdorff on the whole of the
algebras of generalized functions, would by necessity become discrete, therefore trivial, when
restricted to usual, standard smooth functions, see for details Biagioni, Rosinger \& Van der
Walt. \\

Here, in order to further clarify the issue of the possible limitations of the usual
Hausdorff-Kuratowski-Bourbaki concept of topology, let us point out the following. Fundamental
results in Measure Theory, predating the mentioned concept of topology, yet having a clear
topological nature, have never been given a suitable formulation within that
Hausdorff-Kuratowski-Bourbaki concept. Indeed, such is the case, among others, with the
Lebesgue dominated convergence theorem, with the Lusin theorem on the approximation of
measurable functions by continuous ones, and with the Egorov theorem on the relation between
point-wise and uniform convergence of sequences of measurable functions. \\

Similar limitations of the Hausdorff-Kuratowski-Bourbaki concept of topology appeared in the
early 1950s, when attempts were made to turn the convolution of Schwartz distributions into an
operation simultaneously continuous in both its arguments. More generally, it is well known
that, given a locally convex topological vector space, if we consider the natural bilinear
form defined on its Cartesian product with its topological dual, then there will exist a
locally convex topology on this Cartesian product which will make the mentioned bilinear form
simultaneously continuous in both of its variables, if and only if our original locally convex
topology is in fact as particular, as being a normed space topology, see Rosinger \& Van der
Walt for further details. \\

It is also well known that in the theory of ordered spaces, and in particular, ordered groups
or vector spaces, there are important concepts of convergence, completeness, roundedness, etc.,
which have never been given a suitable formulation in terms of the
Hausdorff-Kuratowski-Bourbaki concept of topology. In fact, as seen in Oberguggenberger \&
Rosinger, powerful general results can be obtained about the existence of generalized
solutions for very large classes of nonlinear PDEs, by using alone order structures and their
Dedekind type order completions, without any recourse to any sort of possibly associated
topologies. And the generalized solutions thus obtained can be assimilated with usual
measurable functions, or they can be even more regular, such as being Hausdorff-continuous,
see Rosinger [12]. \\

Finally, it should be pointed out that, recently, differential calculus was given a new
re-foundation by using standard concepts in category theory, such as naturalness. This
approach also leads to topological type processes, among them the so called toponomes or
$\mathcal{C}$-spaces, which prove to be more general than the usual
Hausdorff-Kuratowski-Bourbaki concept of topology, see Nel, and the references cited there. \\

In this way, we can conclude that Mathematics contains a variety of important
\textit{topological type processes} which, so far, could not be formulated in convenient terms
using the Hausdorff-Kuratowski-Bourbaki topological concept. And the differential algebras of
generalized functions, just as much as the far simpler nonstandard reals $^*\mathbb{R}$,
happen to exhibit such a class of topological type processes. \\

On the other hand, the topological type processes on the nowhere dense differential algebras
of generalized functions, used in Mallios \& Rosinger [1], for instance, as well as on the
space-time foam or special space-time foam differential algebras of generalized functions
employed in this paper, see also Mallios [2], Mallios \& Rosinger [2,3], can be given a
suitable formulation, and correspondingly, treatment, by noting that the mentioned algebras
are in fact \textit{reduced powers} see Lo\u{s}, or Bell \& Slomson, of
$\mathcal{C}^{\infty}(X)$, and thus of $\mathcal{C}(X)$ as well. Let us give some further
details related to this claim in the case of the space-time foam algebras. The case of the
nowhere dense algebras was treated in Mallios \& Rosinger [1]. \\

Let us recall, for instance, the definition in (\ref{label1.12}) of
the multi-foam algebras, and note that it obviously leads to the relations

\begin{equation}\label{label3.20}
    B_{L,\, \mathcal{S}}(X) ~=~ (\mathcal{C}^{\infty}(X))^{\Lambda}/
    \mathcal{J}_{L,\, \mathcal{S}}(X) ~\subseteq~ (\mathcal{C}(X))^{\Lambda}/
    \mathcal{J}_{L,\, \mathcal{S}}(X) ~\subseteq~ \mathcal{C}(\Lambda \times X)/
    \mathcal{J}_{L,\, \mathcal{S}}(X)
\end{equation}

assuming in the last term that on $\Lambda$ we consider the discrete topology. A similar
situation holds for the special multi-foam algebras $A_{L,\, \mathcal{S}}(X)$, see
(\ref{label2.14}). \\

Now it is well known, Gillman \& Jerison, that the algebra structure of $\mathcal{C}(\Lambda
\times X)$ is connected to the topological structure of $\Lambda \times X $, however, this
connection is rather sophisticated, as essential aspects of it involve the Stone-\u{C}ech
compactification $\beta(\Lambda \times X)$ of $\Lambda \times X$. And in order to complicate
things, in general $\beta(\Lambda \times X) \neq \beta(\Lambda) \times \beta(X)$, not to
mention that $\beta(\Lambda)$ alone, even in the simplest nontrivial case of $\Lambda =
\mathbb{N}$, has a highly complex structure. \\

It follows that a good deal of the discourse, and in particular, the topological type one, in
the space-time foam and special space-time foam algebras may be captured by the topology of
$\Lambda \times X$, and of course, by the far more involved topology of $\beta(\Lambda \times
X)$. Furthermore, the differential properties of these algebras will, in view of
(\ref{label1.17}) and (\ref{label2.9}), be reducible termwise to classical differentiation of
sequences of smooth functions. \\

In short, in the case of the mentioned differential algebras of generalized functions, owing
to their structure of reduced powers, one obtains a "two-way street" along which, on the one
hand, the definitions and operations are applied to sequences of smooth functions, and then
reduced termwise to such functions, while on the other hand, all that has to be done in a way
which will be compatible with the "reduction" of the "power" by the quotient constructions in
(\ref{label1.12}), or in other words, (\ref{label3.20}) and similarly for (\ref{label2.14}).
By the way, such a "two-way street" approach has ever since the 1950s been fundamental in the
branch of Mathematical Logic, called Model Theory, see Lo\u{s}. But in order not to become
unduly overwhelmed by ideas of Model Theory, let us recall here that the classical
Cauchy-Bolzano construction of the real numbers $\mathbb{R}$ is also a reduced power. Not to
mention that a similar kind of reduced power construction - in fact, its particular case
called "ultra-power" - gives the nonstandard reals $^*\mathbb{R}$ as well. \\ \\

\textbf{Remark 2.}

\medskip
Lately, there has been a growing interest in \textit{noncommutative} studies, and in
particular, noncommutative algebras, see Connes. It is therefore appropriate to mention
possible connections between such noncommutative methods and the space-time foam and special
space-time foam differential algebras of generalized functions in this paper.
\\
In this regard, we recall that, as mentioned at the beginning of this section, in case our
constructions start with arbitrary Banach algebra valued, and not merely real or complex
valued functions, then the resulting space-time foam and special space-time algebras can still
be constructed in the same way, and they will become noncommutative in general. \\

On the other hand, the emergence of noncommutative studies need not at all mean the loss of
interest in, and relevance of commutative structures. Indeed in many problems the commutative
approach turns out to be both more effective and also, of course, much more simple. \\

Finally, it is important to mention here that in the case of singularities of generalized
functions, that is, of singularities in a differential context, the approach in Connes falls
far short even of the long establish linear theory of Schwartz distributions. Indeed, in such
a context, the only differential type operation in Connes, see pp. 19-28, 287-291, is defined
as the commutator with a fixed operator. In this way, it is a rather particular derivation,
even when considered within Banach algebras. The effect is that, it can only to a small extent
deal with the singularities, even when compared with the limited linear Schwartz theory. And
certainly, the approach in Connes can deal with even less with singularities on arbitrary
closed nowhere dense sets, let alone, on the far larger class of arbitrary dense sets whose
complementaries is still dense, such as those in this paper.


\bibliographystyle{plain}

\begin{thebibliography}{99}


\bibitem{}Bell J L, Slomson A B : Models and Ultraproducts, An Introduction.
North-Holland, Amsterdam, 1969

\bibitem{} Biagioni H A : A Nonlinear Theory of Generalized Functions. Lecture
Notes in Mathematics, vol. 1421, Springer, New York, 1990

\bibitem{} Colombeau J-F : New Generalized Functions and Multiplication of Distributions.
Mathematics Studies, vol. 84, North-Holland, Amsterdam, 1984

\bibitem{} Connes A : Noncommutative Geometry. Acad. Press, New York, 1994

\bibitem{} Finkelstein D : Past-future asymmetry of the gravitational field of a
point particle. Physical Review, vol. 110, no. 4, May 1953, 965-967

\bibitem{} Geroch R [1] : What is a singularity in General Relativity ? . Annals of
Physics, vol. 48, 1968, 526-540

\bibitem{} Geroch R [2] : Einstein algebras. Commun. Math. Phys., Springer, vol.
26, 1972, 271-275

\bibitem{} Gillman L, Jerison M : Rings of Continuous Functions. Van Nostrand,
New York, 1960

\bibitem{} Grosser M, Kunzinger M, Oberguggenberger M, Steinbauer R : Geometric Theory of
Generalized Functions with Applications to General Relativity. Kluwer, Dordrecht, 2002

\bibitem{} Gruszczak J, Heller M : Differential Structure of space-time and its
prolongations to singular boundaries. Intern. I. Theor. Physics,
vol. 32, no. 4, 1993, 625-648

\bibitem{} Heller M [1] : Algebraic foundations of the theory of differential spaces.
Demonstratio Math., 24, 1991, 349-364

\bibitem{} Heller M [2] : Einstein algebras and general relativity. Intern. J. Theor.
Physics, vol. 31, no. 2, 1992, 277-288

\bibitem{} Heller M [3] : Thoeretical Foundations of Cosmology, Introduction
to the Global Structure of Space-Time. World Scientific, Singapore,
London, 1092

\bibitem{} Heller M Multarzynski P, Sasin W :  The algebraic approach to spacetime
geometry. Acta Cosmologica, fast. XVI, 1989, 53-85

\bibitem{} Heller M, Sasin W [1] : Generalized Friedman's equation and its singularities.
Acta Cosmologica, fasc. XIX, 1993: 23-33

\bibitem{} Heller M, Sasin W [2] : Sheaves of Einstein algebras. Intern. J. Theor.
Physics, vol. 34, no. 3, 1995, 387-398

\bibitem{} Heller M, Sasin W [3] : Structured spaces and their application to
relativistic physics. J. Math. Phys., 36, 1995, 3644-3662

\bibitem{} Kaneko A : Introduction to Hyperfunctions. Kluwer, Dordrecht, 1088

\bibitem{} Kirillov A A [1] : Elements of the Theory of Representations. Springer,
New York, 1976

\bibitem{} Kirillov A A [2] : Geometric Quantization. In (Eds. Arnold V I,Novikov
S P) Dynamical Systems IV. Symplectic Geometry and its Application.
Springer, New York, 1990, 137-172

\bibitem{} Lo\u{s} J : On the categoricity in power of elementary deductive systems
and some related problems. Colloq. Math., 3, 1954, 58-62

\bibitem{} Mallios A [1] : Geometry of Vector Sheaves. An Axiomatic Approach
to Differential Geometry, vols. I (chapts. 1-5), II (chapts. 6-11)
Kluwer, Amsterdam: 1998

\bibitem{} Mallios A [2] : Modern Differential Geometry in Gauge Theories. Volume 1 :
Maxwell Fields, Volume 2 : Yang-Mills Fields. Birkhauser, Boston, 2006

\bibitem{} Mallios A [3] : On an axiomatic treatment of differential geometry
via vector.sheaves. Applications. (International Plaza) Math.
Japonica, vol. 48, 1998, pp. 93-184

\bibitem{} Mallios A [4] : The de Rham-K\"{a}hler complex of the Gel'fand sheaf of
a topological algebra. J. Math. Anal. Appl., 175, 1993, 143-168

\bibitem{} Mallios A [5] : On an abstract form of Weil's intergrality theorem. Note
Mat., 12, 1992, 167-202 (invited paper)

\bibitem{} Mallios A [6] : On an axiomatic approach to geometric prequantization
: A classification scheme \'a la Kostant-Souriau-Kirillov. J.
Mathematical Sciences (former J. Soviet Math.), vol. 98-99

\bibitem{} Mallios A, Rosinger E E [1] : Abstract differential geometry, differential
algebras of generalized functions, and de Rham cohomology. Acta
Applicandae Mathematicae (accepted)

\bibitem{} Mallios A, Rosinger E E [2] : Space-time foam dense singularities and
de Rham cohomology   (to appear)

\bibitem{} Mallios A, Rosinger E E [3] : Dense singularities and de Rham cohomology. In (Eds.
Strantzalos P, Fragoulopoulou M) Topological Algebras with Applications to Differential
Geometry and Mathematical Physics. Proc. Fest-Colloq. in honour of Prof. Anastasios Mallios
(16-18 September 1999), pp. 54-71, Dept. Math. Univ. Athens Publishers, Athens, Greece, 2002

\bibitem{} Mostow M A : The differentiable space structures of Milnor classifying
spaces, simplicial complexes, and geometric realizations. J. Diff.
Geom., 14, 1979, 255-293

\bibitem{} Nel L D : Differential calculus founded on an isomorphism. Appl. Categorial
Structures, 1, 1993, 51-57

\bibitem{} Oberguggenberger M B : Multiplication of Distributions and Applications to PDEs.
Pitman Research Notes in Mathematics, Vol. 259, Longman, Harlow, 1992

\bibitem{} Oberguggenberger M B, Rosinger E E : Solution of Continuous Nonlinear
PDEs through Order Completion. Mathematics Studies, vol. 181,
North-Holland, Amsterdam, 1994

See also review MR 95k:35002

\bibitem{} Oxtoby J C : Measure and Category. Springer, New York, 1971

\bibitem{} Rosinger E E [1] : Embedding of the $\mathcal{D}'$ distributions into
pseudotopological algebras. Stud. Cerc. Mat., vol. 18, no. 5, 1966, 687-729

\bibitem{} Rosinger E E [2] :  Pseudotopological spaces, the embedding of the
$\mathcal{D}'$ distributions into algebras. Stud. Cerc. Mat., vol. 20, no. 4, 1968, 553-582

\bibitem{} Rosinger E E [3] : Distributions and Nonlinear Partial Differential Equations.
Lecture Notes in Mathematics, vol. 684, Springer: New York, 1978

\bibitem{} Rosinger E E [4] : Nonlinear Partial Differential Equations, Sequential
and Weak Solutions. Mathematics Studies, vol. 44, North-Holland, Amsterdam, 1980

\bibitem{} Rosinger E E [5] : Generalized Solutions of Nonlinear Partial Differential
Equations. Mathematics Studies, vol. 146, North-Holland, Amsterdam, 1987

\bibitem{} Rosinger E E [6] : Nonlinear Partial Differential Equations, An A1gebraic View of
Geueralized solutions. Mathematics Studies, vol. 164, North-Holland, Amsterdam, 1990

\bibitem{} Rosinger E E [7] : Global Version of the Cauchy-Kovalevskaia Theorem for Nonlinear
PDEs. Acta Applicandae Mathematicae, Vol. 21, 1990, pp. 331-343

\bibitem{} Rosinger E E [8] : Parametric Lie Group Actions on Global Generalized Solutions of
Nonlinear PDEs, including a Solution to Hilbert's Fifth Problem, Kluwer Acad. Publ.,
Dordrecht, Boston, London, 1998

\bibitem{} Rosinger E E [9] : Space-time foam differential algebras of generalized
functions. Private communication. Vancouver, 1998

\bibitem{} Rosinger E E [10] : Dense Singularities and Nonlinear PDEs (to appear)

\bibitem{} Rosinger E E [11] : Differential algebras with dense singularities on manifolds.
Acta Applicandae Mathematicae, Vol. 95, No. 3, Feb. 2007, 233-256, arXiv:math.DG/0606358

\bibitem{} Rosinger E E [12] : Can there be a general nonlinear PDE theory for the existence
of solutions ? math.AP/0407026

\bibitem{} Rosinger E E [13] : Singularities and flabby sheaves. (to appear)

\bibitem{} Rosinger E E [14] : Scattering in highly singular potentials. \\
arXiv:quant-ph/0405172

\bibitem{} Rosinger E E [15] : Which are the Maximal Ideals ? \\ arXiv:math.GM/0607082

\bibitem{} Rosinger E E, Van der Walt J-H : Beyond topologies ( to appear)

\bibitem{} Rosinger E E, Walus Y E [1] : Group invariance of generalized solutions obtained
through the algebraic method. Nonlinearity, Vol. 7, 1994, pp. 837-859

\bibitem{} Rosinger E E, Walus Y E [2] : Group invariance of global generalized solutions of
nonlinear PDEs in nowhere dense algebras. Lie Groups and their Applications, Vol. 1, No. 1,
July-August 1994, pp. 216-225

See also reviews : MR 92d:46008, Zbl. Math. 717 35001, MR 92d:46097,
Bull. AMS vo1.20, no.1, Jan 1989, 96-101, MR 89g:35001

\bibitem{44} Sasin W [1] : The de Rham cohomology of differential spaces. Demonstration
Mathematica, vol. XXII, no. 1, 1989, 249-270

\bibitem{45} Sasin W [2]   : Differential spaces and singularities in differential spacetime.
Demonstratio Mathematica, vol. XXIV, no. 3-4, 1991, 601-634

\bibitem{46} Sikorski R : Introduction to Differential Geometry (in Polish). Polish
Scientific Publishers. Warsaw, 1972

\bibitem{47} Souriau J-M [1] : Structures des Syst\`{e}mes Dynamiques. Dunod, Paris,
1970

\bibitem{48} Souriau I-M [2] : Groupes Diff\'erentiels. In Differential Geometric
Methods in Mathematical Physics. Lecture Notes in Mathematics vol.
836, Springer, New York, 1980, 91-128

\bibitem{} Synowiec J A : Some highlights in the development of Algebraic Analysis.
Algebraic Analysis and Related Topics. Banach Center Publications, Vol. 53, 2000, 11-46,
Polish Academy of Sciences, Warszawa

\end{thebibliography}

\end{document}